\documentclass{article}

\usepackage[a4paper,hmargin=3cm,vmargin=3cm]{geometry}

\usepackage{amsmath}
\usepackage{amsfonts}
\usepackage{amssymb}
\usepackage{epsfig}
\usepackage{theorem}
\usepackage{comment}
\usepackage[dvips]{color}
\usepackage{color}
\usepackage{graphicx}

\def\be{\begin{equation}}
\def\ee{\end{equation}}

\def\C{{\mathbb C}} 
\def\N{{\mathbb N}} 

\def\Z{{\mathbb Z}}
\def\R{{\mathbb R}} 
\def\Q{{\mathbb Q}}

\def\phi{{\varphi}}
\def\v{{\varepsilon}} 
\def\deg{{\rm deg\,}}

\def\cos{{\rm cos\,}} 
 
\def\GCD{{\rm GCD }}

\def\qed{$\ \ \Box$ \vskip 0.2cm}
\def\t{\widetilde}

\def\tilde{\widetilde}

\def\bp{\begin{proposition}}
\def\ep{\end{proposition}}

\def\bt{\begin{theorem}}
\def\et{\end{theorem}}
\def\br{\begin{remark}}
\def\er{\end{remark}}
\def\be{\begin{equation}}
\def\bee{\begin{equation*}}
\def\la{\label}
\def\l{\label}

\def\ee{\end{equation}}
\def\eee{\end{equation*}}
\def\bl{\begin{lemma}}
\def\el{\end{lemma}}
\def\br{\begin{remark}}
\def\er{\end{remark}}
\def\bc{\begin{corollary}}
\def\ec{\end{corollary}}
\def\pr{\noindent{\it Proof. }}

\def\bd{\begin{definition}}
\def\ed{\end{definition}}
\input epsf.sty

\newtheorem{theorem}{Theorem}[section]
\newtheorem{lemma}[theorem]{Lemma}
\newtheorem{corollary}[theorem]{Corollary}
\newtheorem{proposition}[theorem]{Proposition}
\newtheorem{remark}{Remark}[section]

\begin{document}

\title{Solution of the parametric center problem for the Abel differential equation} 
\author{F. Pakovich}
\date{}

\maketitle
\begin{abstract}
The Abel differential equation $y'=p(x)y^2+q(x)y^3$
with $p,q\in \R[x]$ 
is said to have a center on a segment $[a,b]$ if all its solutions, with the initial
value $y(a)$ small enough, satisfy the condition $y(b)=y(a)$. The problem of description of conditions implying that the Abel equation has a center may be interpreted as a simplified version of the classical Center-Focus problem of Poincar\'e.
The Abel equation is said to have a ``parametric center" if for each $\v \in \R$ the equation 
$y'=p(x)y^2+\v q(x)y^3$
has a center. In this paper we show that the Abel equation has a parametric center if and only if the
antiderivatives 
$P=\int p(x) dx,$ $Q=\int q(x) dx$ satisfy the equalities  
$P=\t P \circ W,\ $ $Q=\t Q\circ W$ for 
some polynomials $\t P,$ $\t Q,$ and $W$ such that $W(a)=W(b)$.
We also show that the last condition is necessary and sufficient for the ``generalized moments''
$\int_a^b P^id Q$ and   
$\int_a^b Q^id P$ to vanish for all $i\geq 0.$

\end{abstract}

\section{Introduction}
Let 
\be \la{abc} y'=p(x)y^2+q(x)y^3\ee
be the Abel differential equation, where $x$ is real and $p(x)$ and $q(x)$ are continuous. 
Equation \eqref{abc}
is said to have a center on a segment $[a,b]$ if all its solutions, with the initial
value $y(a)$ small enough, satisfy the condition $y(b)=y(a)$. 

The problem of description of conditions implying a center for \eqref{abc}  is closely related with 
the classical Poincar\'e center-focus problem about conditions implying that all solutions of the system 
\be \la{sys1}
\left\{
\begin{array}{rcl}
\dot{x} & =& -y+F(x,y), \\ 
\dot{y} & =&x+G(x,y),\\
\end{array} 
\right.
\ee
where $F(x,y),$ $G(x,y)$ are polynomials without constant and linear terms,
around zero are closed. Namely, it was shown in \cite{cher} that if $F(x,y),$ $G(x,y)$ are homogeneous polynomials of the same degree, then 
one can construct trigonometric polynomials $f(\cos \phi,\sin \phi),$ $g(\cos \phi,\sin \phi)$ such that \eqref{sys1} has a center if and only if all
solutions of the equation 
$$
\frac{d r}{d \phi}=f(\cos \phi,\sin \phi)\, r^2+g(\cos \phi,\sin \phi)\, 
r^3 
$$
with $r(0)$ small enough
are periodic on $[0,2 \pi]$.

Set
\be \label{form} P(x)=\int_0^{x}p(s)ds , \ \ \ Q(x)=\int_0^{x} q(s)ds.\ee 
The following ``composition condition'' introduced in \cite{al} is  sufficient 
for equation \eqref{abc} to have a center:
there exist $C^1$-functions $\tilde P,\tilde Q, W$ such that  
\be \label{c} P(x)=\tilde P(W(x)), \ \ \ \ Q(x)=\tilde Q(W(x)),\ \ \ \ W(a)=W(b). \ee
Indeed, if \eqref{c} holds, then any solution of \eqref{abc} has the form $y(x)=\tilde y(W(x)),$ where 
$\tilde y$ is a solution of the equation 
$$ y'=\tilde P^{\prime}(x)y^3+\tilde Q^{\prime}(x)y^2,$$
implying that $y(a)=y(b)$, since $W(a)=W(b).$
\pagebreak

It is known that in general the composition condition is not necessary for \eqref{abc} to have a center (\cite{alw}). However, 
it 
is believed that, in the case where $p(x)$ and $q(x)$ are polynomials, equation \eqref{abc} has a center if and only if the 
composition condition \eqref{c}  
holds for some polynomials $\t P, \t Q, W\in \R[x]$
(see  \cite{bpy}, \cite{bry} for some  partial results in this direction). 

In this paper we study the following ``parametric center problem''
for equation \eqref{abc} 
with polynomial coefficients: 
{\it under what conditions the equation
\be \la{abp} y'=p(x)y^2+\v q(x)y^3, \ \ \  p,q\in \R[x],
\ee  
has a center for any $\v \in \R$?} Posed for the first time in the series of papers \cite{bfy2}, \cite{bfy3}, \cite{bfy4},  this problem turned out to be very constructive
and resulted in a whole area of new ideas and methods related to the so called ``polynomial moment problem'' 
(see the discussion below). 
However, in its full generality the parametric center problem  remained unsolved (see the recent paper \cite{bpy1} for the state of the art), and the goal of this paper is to fill this gap. 
Our main result is the following theorem.

\bt \la{t1} Abel differential equation \eqref{abp}  
has a center on a segment $[a,b]$ for any $\v \in \R$ if and only if the antiderivatives $P=\int p(x) dx$ and $Q=\int q(x) dx$ 
satisfy  composition condition \eqref{c} for some polynomials $\t P,$ $\t Q$, $W$.
\et

The proof of Theorem \ref{t1} is based on the link between the parametric center problem and the ``polynomial moment problem''. Namely, it was shown in \cite{bfy4} that the parametric center implies the equalities
\be \la{mix} 
\int_a^b P^id Q=0, \ \ \ i\geq 0, \ \ \  \int_a^b Q^id P=0, \ \ \ i \geq 0.
\ee
We will call the problem of description of solutions \eqref{mix} the ``mixed polynomial moment problem'', and the problem of 
description of solutions 
\be \la{mom} 
\int_a^b P^id Q=0,  \ \ \ i\geq 0,
\ee
simply the ``polynomial moment problem''.

The polynomial moment problem 
has been studied in many recent papers (see e.g. \cite{bfy2},\cite{bfy3} \cite{bfy4}, \cite{c}, \cite{acoun}, 
\cite{pa2}, \cite{pp}, \cite{pry}, \cite{mp1}, \cite{pak1}, 
\cite{pak}, \cite{ro}).  Again, the composition condition \eqref{c} is sufficient for equalities \eqref{mom} to be satisfied although in general is not necessary (\cite{acoun}). A complete solution of the polynomial moment problem was obtained in the recent papers \cite{mp1}, \cite{pak}. Namely, it was shown in \cite{mp1} that if polynomials $P,$  $Q$ satisfy \eqref{mom}, then there exist polynomials $Q_j$ such that $Q=\sum_j Q_j$ and
\be \l{cc}
P(x)=P_j(W_j(x)), \ \ \
Q_j(x)=V_j(W_j(x)), \ \ \ W_j(a)=W_j(b)
\ee
for some polynomials $P_j(z), V_j(z), W_j(z)$. Moreover, in 
\cite{pak} polynomial solutions of \eqref{mom} were desc\-ribed in 
an explicit form  (see Section 2 below).

In this paper we apply results of \cite{pak} to each of the two systems appeared in \eqref{mix} separately and show that the restrictions obtained 
imply that any solution $P,Q$ of the mixed polynomial moment problem satisfy  composition condition \eqref{c}. Thus, in fact we prove the following ``moment'' 
counterpart of Theorem \ref{t1}.

\bt \la{t2} Polynomials $P, Q\in \R[x]$  satisfy equalities \eqref{mix} if and only if they 
satisfy composition condition \eqref{c} for some polynomials $\t P, \t Q, W\in \R[x].$
\et

Although the center problem for the Abel equation with polynomial coefficients can be consi\-dered in the complex setting, in this paper we  
work in the classical real framework. Thus,  
we start the paper from the adaptation for the real case of the results of \cite{pak}, obtained over $\C$. 
Namely, we show in Section 2 that possible ``types'' of solutions of the polynomial moment problem over $\R$ remain the same, although one of these types becomes ``smaller'' (Theorem \ref{mpr}). 
Besides, in Section 2 we establish some important restrictions of the arithmetical nature on points $a,b$ for which there exist solutions of \eqref{mom} 
which do not satisfy 
the composition condition  (Corollary \ref{red}).

In Section 3 we apply  the results of Section 2 to \eqref{mix}, and prove Theorem \ref{t2}.
The main difficulties of the proof 
stem from the fact that after solving systems in \eqref{mix} separately we arrive to 
functional equations of the type 
\be \la{eba} \sum_{j=1}^r V_j(W_j(x))=A(B(z)),\ee where $V_j,W_j,A,B$ are polynomials, and $r$ equals $2$ or 3. Such equations can be considered as generalizations of the functional equation  
\be \la{ebaq} A(B(x))=C(D(x)),\ee studied by Ritt (\cite{rit}).
However, the well established methods for studying \eqref{ebaq}, related to the monodromy, cannot be applied to \eqref{eba} for $r>1$, and essentially the only method which remains is a  painstaking analysis of coefficients. Such an analysis in general leads to rather
complicated systems of equations, and Theorem \ref{t2} is deduced from restrictions on $P$ and $Q$ obtained from these systems combined with restrictions on possible values of $a$ and $b$.

\section{Polynomial moment problem over $\C$ and over $\R$}
\subsection{Solution of the polynomial moment problem over $\C$}
In this subsection we briefly recall a description of
$P,Q\in \C[z]$  satisfying \eqref{mom} for $a,b\in \C$. For more details we refer 
the reader to \cite{pak}.

Recall that the Chebyshev polynomials (of the first kind) $T_n$ are defined by the formula $T_n(\cos \phi)=\cos(n \phi).$ It follows directly from this definition 
that 
\be \la{tyui} T_n(1)=1, \ \ \ \ T_n(-1)=(-1)^n, \ \ \ n\geq 0\ee
and 
$$T_n\circ T_m=T_m\circ T_n=T_{mn},\ \ \ n,m\geq 1, $$
where the symbol $\circ$ denotes a composition of functions, $A\circ B=A(B(z))$.

An explicit expression for $T_n$ is   
given by the formula
\be \la{cheb} T_n=\frac{n}{2}\sum_{k=0}^{[n/2]}(-1)^k\frac{(n-k-1)!}{k!(n-2k)!}(2x)^{n-2k},\ee 
implying in particular that 
\be \la{kots} T_n(-x)=(-1)^nT_n(x).\ee

Following \cite{pak}, we will call a solution $P, Q$ of \eqref{mom} reducible if composition condition \eqref{c} holds 
for some $\t P, \t Q, W\in \C[x].$

\bt \l{mpc} (\cite{pak}) Let $P,$ $Q$ be non-constant complex polynomials
and $a,b$ distinct complex numbers such that equalities \eqref{mom} hold. Then, either $Q$ is a reducible solution of \eqref{mom}, or there exist 
complex polynomials $P_j,$ $Q_j,$ $V_j,$ $W_j,$ $1\leq j \leq r,$ 
such that 
$$Q=\sum_{j=1}^rQ_j, \ \ \ 
P=P_j\circ W_j, \ \ \
Q_j=V_j\circ W_j, \ \ \  \ \ \ W_j(a)=W_j(b).
$$ 
Moreover,
one of the following conditions holds:
\vskip 0.25cm 
\noindent 1) $r=2$ and 
$$P=U\circ z^{sn}R^n(z^n) \circ V,\ \ \ W_1=z^n \circ V, \ \ \ W_2=z^sR(z^n)\circ V,$$ 
where $R,$ $U$, $V$  are complex polynomials, $n>1$, $s> 0,$ $\GCD(s,n)=1;$

\vskip 0.25cm 
\noindent 2)  $r=2$ and
$$P=U\circ T_{m_1m_2} \circ V,\ \ \ W_1= T_{m_1} \circ V, \ \ \  
W_2= T_{m_2}\circ V,$$ 
where $U$, $V$ are complex polynomials, $m_1>1,$ $m_2>1$, $\GCD(m_1,m_2)=1;$

\pagebreak

\vskip 0.25cm 
\noindent 3) $r=3$ and 
$$P=U\circ z^{2}R^2(z^{2})\circ  T_{m_1m_2}
\circ V,$$ 
$$W_1=T_{2m_1} \circ V, \ \ \ W_2=T_{2m_2}\circ V, \ \ \ 
W_3= (zR(z^2) \circ T_{m_1m_2}) \circ V,$$ 
where $R,$ $U$, $V$ are complex polynomials, 
$m_1>1,$ $m_2>1$ are odd, and $\GCD(m_1,m_2)=1.$
\et

We will call solutions 1), 2), 3) appearing in Theorem \ref{mpc} solutions of the first, the second, and the third type 
correspondingly.

Notice that these sets of solutions are not disjointed. For example, if one of the parameters $n$, $m$ of a solution of the second type equals 2, then this solution is also 
a solution of the first type. Indeed, if say $n=2,$ then $W_1=T_2\circ V=\mu\circ z^2\circ V,$ where $\mu =2x-1.$ On the other hand,  since $m$ is odd in view 
of $\GCD(n,m)=1,$ the polynomial $W_2=T_m\circ V$  has the form $T_m=zR(z^2)\circ V$ by \eqref{cheb}. 
Therefore, 
$$P=(U\circ \mu)\circ z^{2}R^2(z^{2})\circ V,\ \ \  Q=((V_1 \circ \mu)\circ z^2+V_2\circ z R(z^2))\circ V,$$
and for the polynomials $\t W_1=z^2\circ V$ and $\t W_2=W_2=z R(z^2)\circ V$ the equalities 
\be \la{aqwer} \t W_1(a)=\t W_1(b), \ \ \ \ \t W_2(a)=\t W_2(b)\ee hold.

Similarly, if the 
parameters $a,b$ of a solution of the third type satisfy  $V(a)=-V(b)$, then  this solution is also a 
solution of the first type. Indeed, $$V_1\circ T_{2m_1}+V_2\circ T_{2m_2}=\t V_1 \circ z^2$$ for some $\t V_1\in \C[z]$, while $$zR(z^2) \circ T_{m_1m_2}= z\t R(z^2)$$ for some $\t R\in \C[z]$, since $m_1,$ $m_2$  are odd. Therefore, 
$$P=U\circ z^{2}\t R^2(z^{2})\circ V,\ \  \ Q=(\t V_1 \circ z^2+V_3\circ z\t R(z^2))\circ V,$$ and $\t W_1=z^2\circ V$ satisfies 
$\t W_1(a)= \t W_1(b)$, since $V(a)=-V(b)$.

Finally, one can check that a solution of the third type is a solution of the second type if 
$(T_{m_i}\circ V)(a)=(T_{m_i}\circ V)(b)$ for $i$ equals 1 or 2 (for more details concerning interrelations between different types of solutions see \cite{pak}, pp. 725-726).

\subsection{Lemmas related to $a,b$}
In this subsection we collect some results implying restrictions on the integration limits  $a,b$ appearing in solutions of the second and the third types.

\bl \la{skun} Let $T_{m_1},$ $T_{m_2},$ $T_{m_3}$ be the Chebyshev polynomials and $a$, $b$ be distinct complex numbers. 
\vskip 0.1cm
\noindent a) Assume that 
\be \la{volk} T_{m_1}(a)=T_{m_1}(b), \ \ \ \ \ \ T_{m_2}(a)=T_{m_2}(b).\ee 
Then either  $T_l(a)=T_l(b)$ for $l=\GCD(m_{1},m_{2}),$ or \be \la{uraa} T_{m_1m_2}^{\prime }(a)=T_{m_1m_2}^{\prime }(b)=0.\ee
\vskip 0.1cm
\noindent b) Assume that \be \la{uuiioo} T_{m_1}(a)=T_{m_1}(b), \ \ \ \ \ \ T_{m_2}(a)=T_{m_2}(b), \ \ \ \ \ \ T_{m_3}(a)=T_{m_3}(b).\ee
Then there exists a pair of distinct indices $i_1,i_2,$ $1\leq i_1,i_2 \leq 3,$ 
such that  $T_l(a)=T_l(b)$ for $l=\GCD(m_{i_1},m_{i_2})$.

\el

\pr Choose $\alpha, \beta \in \C$ such that $\cos \alpha= a,$
$\cos \beta =b.$ Then equalities \eqref{volk}
imply the equalities
\be \la{100} m_1 \alpha = \v_1 m_1 \beta + 2\pi k_1, \ \ \
m_2 \alpha = \v_2 m_2 \beta+ 2\pi k_2, \ee where $\v_1=\pm 1,$ $\v_2=\pm 1,$ 
and $k_1,k_2\in \Z.$ Assume first that $\v_1=\v_2.$ Let  $u,v$ be integers satisfying \be\la{uv} um_1+vm_2=l.\ee Multiplying the 
first equality in \eqref{100} by $u$ and adding 
the second equality multiplied by $v$, 
we see that $$l\alpha= \v_1l\beta + 2\pi k_1 u+2\pi k_2 v,$$  implying that 
$T_{l}(a)=T_{l}(b)$.

Assume now that  $\v_1=-\v_2.$ Then, similarly, we conclude that 
\be \la{qaz} l\alpha= \v_1\beta(um_1-vm_2) + 2\pi k_1 u+2\pi k_2 v.\ee
Furthermore, eliminating $\alpha$ from equalities \eqref{100} we obtain \be \la{vu} \v_1m_1m_2\beta=\pi k_2 m_1-\pi k_1 m_2.\ee Since 
\be \la{pizz} T_n^{\prime}(\cos \phi)=n(\sin n\phi/\sin \phi),\ee equality \eqref{vu} implies that $T_{m_1m_2}^{\prime }(b)=0$, unless 
$\beta=\pi k_3,$ $k_3\in \Z.$ In the last case, $b=1$ if $k_3$ is even, and $b=-1$ if $k_3$ is odd, implying that
$$T_{l}(b)=(-1)^{k_3l},$$ in view of \eqref{tyui}. On the other hand, if $\beta=\pi k_3$, then
\eqref{qaz} implies that  
$$T_{l}(a)=(-1)^{k_3(u m_1-v m_2)}.$$ Since the sum and the difference of any two numbers have the same parity this implies that 
$$T_{l}(a)=(-1)^{k_3(u m_1+v m_2)}=(-1)^{k_3l}=T_{l}(b).$$ Similarly, one can see that 
$T_l(a)=T_l(b)$ unless $T_{m_1m_2}^{\prime }(a)=0.$

In order to prove b) observe that equalities \eqref{uuiioo} imply the equalities   
\be \la{10} m_1 \alpha = \v_1 m_1 \beta + 2\pi k_1, \ \ \
m_2 \alpha = \v_2 m_2 \beta+ 2\pi k_2, \ \ \ 
m_3 \alpha = \v_3  m_3 \beta+ 2\pi k_3, \ee where $\v_1=\pm 1,$ $\v_2=\pm 1,$ $\v_3=\pm 1,$ 
and $k_1,k_2,k_3\in \Z.$ 
Clearly, among the numbers $\v_1,$ $\v_2,$ $\v_3$ at least two are equal 
and we conclude as above that the equality $T_{l}(a)=T_{l}(b)$ holds for $l=\GCD(m_{i_1},m_{i_2})$, where $\v_{i_1}=\v_{i_2}.$
\qed

\bc \la{xyi}  Let $T_{m_1},$ $T_{m_2}$ be the Chebyshev polynomials  and $a$, $b$ be distinct complex numbers. 
Assume that 
$$T_{m_1}(a)=T_{m_1}(b), \ \ \ \ \ \ T_{m_2}(a)=T_{m_2}(b)$$ and $\GCD(m_1,m_2)=1$.
Then $$ T_{m_1m_2}^{\prime }(a)=T_{m_1m_2}^{\prime }(b)=0.$$
\ec
\pr Follows from Lemma \ref{skun}, a) taking into account that 
the equality $S(a)=S(b)$ for some polynomial $S$ and $a\neq b$ obviously implies that $\deg S>1$. \qed

Recall that a number $\gamma\in \C$ is called algebraic if it is a root of an equation with rational coefficients. The set of all algebraic numbers is a subfield of $\C.$ A monic polynomial $p(x)\in \Q[x]$ of the minimal degree 
such that $p(\gamma)=0$ is called a minimal polynomial of $\gamma.$ A minimal polynomial is irreducible over $\Q$. An algebraic number $\gamma$ is called an algebraic integer  if its minimal polynomial has integer coefficients. In fact, this condition may be replaced by a weaker condition 
that $\gamma$ is a root of {\it  some}  monic polynomial with integer coefficients.
The set of all algebraic integers is closed under addition and multiplication.

\bl \la{poi} Assume that $a\in \C$ is a root of $T_n^{\prime}.$ Then $a\in \R$, and $2a$ is an algebraic integer. 
\el 
\pr Since equality \eqref{pizz} shows that $T^{\prime}_n$ has $n-1$ distinct real roots, all roots of 
$T^{\prime}_n$ are real. The other statements follow 
from the formulas $$T_n^{\prime}=nU_{n-1}$$ and 
$$U_n=\sum_{k=0}^{[n/2]}(-1)^k \binom{n-k}{k}(2x)^{n-2k},$$ 
where $U_n$ denotes the Chebyshev polynomial of the second kind. \qed

\bc \la{red} In the notation of Theorem \ref{mpc} assume that $Q$ is a solution of \eqref{mom} of the second type, or   
a solution of the third type,
which cannot be represented as a solution of the first type. Then $2V(a)$ and $2V(b)$ are algebraic integers. 
\ec

\pr Without loss of generality we may assume that $V=x.$ If $Q$ is of the second type, then the statement follows from Corollary \ref{xyi} and Lemma \ref{poi}.

If $Q$ is of the third type, then applying  
Lemma \ref{skun}, a) and Lemma \ref{poi}
to the equalities 
$$T_{2m_1}(a)=T_{2m_1}(b), \ \ \ \ T_{2m_2}(a)=T_{2m_2}(b),$$ we conclude that $2V(a)$ and 
$2V(b)$ are algebraic integers, unless  $T_2(a)=T_2(b)$. However, 
the last equality yields the equality $a=-b,$ implying, as it was observed above, that  
$Q$ can be represented as a solution of the first type. \qed

\subsection{Decompositions of polynomials with real coefficients}
In this subsection we collect necessary results concerning decomposition of polynomials with real coefficients into compositions of polynomials of lesser degree.

The following lemma is well known (see e. g. Corollary 2.2 in \cite{pak}).

\bl \la{wsx} Assume that $$P=A\circ B=\t A\circ \t B,$$ where $P, A,B,\t A,\t B\in \C[z]$ and
$\deg A=\deg \t A.$ Then there exists a polynomial $\mu\in \C[z]$ of degree one such that 
$$\t A=A\circ \mu^{-1}, \ \ \t B=\mu\circ B.\ \ \ \Box$$
\el

\bc \la{reco0} Let $P=U\circ V,$ where $P\in \R[x],$ while $U,V\in \C[z].$ 
Assume that the leading coefficient of $V$ and its constant term are real numbers. Then 
$U,V\in \R[x].$
\ec
\pr Since $P\in \R[x],$ we have:
\be \la{uip} P=U\circ V=\overline{U}\circ \overline{V},\ee where $\overline{U},$ $\overline{V}$ are polynomials obtained from $U,V$ by the complex conjugation of all coefficients. By Lemma \ref{wsx},  equality \eqref{uip} implies that 
\be \la{iui} \overline{U}=U\circ \mu^{-1}, \ \ \overline{V}=\mu\circ V,\ee
where $\mu=\alpha z+\beta$ for some $\alpha, \beta\in \C.$ Since the leading coefficient of $V$ is real, the second equality in \eqref{iui}
implies that $\alpha=1$. Now the equality 
$\overline{V}=V+\beta$ implies that $\beta=0,$ since the constant term of $V$ is real. 
Therefore, $\overline{U}=U,$ $\overline{V}=V$ and 
hence $U,V\in \R[x].$ \qed

\bc \la{reco} Assume that $$P=U\circ V,$$ where $P\in \R[x],$ while $U,V\in \C[z].$ Then there exists a polynomial $\mu\in \C[z]$ of degree one such that the polynomials 
$$U_1=U\circ \mu^{-1}, \ \ V_1=\mu\circ V$$
are contained in $\R[x].$
\ec
\pr Let $\mu$ be any polynomial of degree one such that the leading coefficient and the constant term
of
the polynomial $V_1=\mu\circ V$ are real numbers.
Then $U_1$ and $V_1$ are contained in $\R[x]$ by Corollary \ref{reco0}. \qed

\bl \la{eche} Let $\mu_1,\mu_2$ be complex polynomials of degree one.
\vskip 0.1cm
\noindent a) Assume that the polynomial $\mu_1\circ z^n\circ \mu_2$, $n\geq 2,$ has real coefficients.
Then there exist $\t \mu\in \R[x]$ and $c\in \C$ such that 
$\mu_2=c\t \mu.$
\vskip 0.1cm
\noindent b) Assume that the polynomial $\mu_1\circ T_n\circ \mu_2$, $n\geq 2,$ has real coefficients.
Then either $\mu_2\in \R[x]$ or 
there exists $\t \mu\in \R[x]$ such that 
$\mu_2=i\t \mu.$ 
\el
\pr 
Let $\mu_1=\alpha_1z+\beta_1,$ $\mu_2=\alpha_2z+\beta_2,$ where $\alpha_1,\beta_1,\alpha_2,\beta_2\in \C.$ Then the coefficients of
$z^n$ and $z^{n-1}$ of the polynomial $\mu_1\circ z^n\circ \mu_2$ are $c_n=\alpha_1\alpha_2^n$ and 
$c_{n-1}=\alpha_1\alpha_2^{n-1}\beta_2 n$, correspondingly. Since by assumption these numbers are real, we conclude that $$ \frac{c_{n-1}}{c_n}=\frac{n\beta_2}{\alpha_2}$$ is also a real number. Therefore, 
$\beta_2/\alpha_2\in \R$, and hence 
$\mu_2=\alpha_2\t \mu,$ where $\t \mu=x+(\beta_2/\alpha_2)\in \R[x].$

Similarly, since 
$$T_n=2^{n-1}x^n-n2^{n-3}x^{n-2}+\dots$$ by \eqref{cheb}, the coefficients of 
$z^n$, $z^{n-1}$, $z^{n-2}$ of the polynomial $\mu_1\circ T_n\circ \mu_2$ are: $$c_n=\alpha_12^{n-1}\alpha_2^n, \ \ \  c_{n-1}=\alpha_12^{n-1}n\alpha_2^{n-1}\beta_2, \ \ \ c_{n-2}=\alpha_12^{n-2}n(n-1)\alpha_2^{n-2}\beta_2^2-\alpha_12^{n-3}n\alpha_2^{n-2},
$$
correspondingly. As above, $c_n,c_{n-1}\in \R$ implies that $\beta_2/\alpha_2\in \R.$ Since
$$c_{n-2}=\frac{n(n-1)c_n}{2}\left(\frac{\beta_2}{\alpha_2}\right)^2-\frac{nc_n}{4}\left(\frac{1}{\alpha_2}\right)^2,$$
it follows now from $c_{n-2}\in \R$ that $\alpha_2^2\in \R$ implying the statement. \qed

\subsection{Solution of the polynomial moment problem over $\R$}
In this subsection we deduce from Theorem \ref{mpc} a description of polynomials
$P,$ $Q$ with {\it real} coefficients satisfying \eqref{mom} for  $a,b\in \R$.

The theorem below is a ``real'' analogue of Theorem \ref{mpc}.
Keeping the above notation we will call a solution $Q\in \R[x]$ of \eqref{mom} reducible if \eqref{c} holds 
for some $\t P, \t Q, W\in \R[x].$ We also will call solutions 1), 2), 3) 
described below solutions of the first, the second, and the third type.
Notice that the set of
solutions of the first type in the real case is ``smaller'' than the one in the complex case.

\bt \l{mpr} Let $P,$ $Q$ be non-constant real polynomials
and $a,b$ distinct real numbers such that equalities \eqref{mom} hold. Then, either $Q$ is a reducible solution of \eqref{mom}, or there exist 
real polynomials $P_j,$ $Q_j,$ $V_j,$ $W_j,$ $1\leq j \leq r,$ 
such that 
$$Q=\sum_{j=1}^rQ_j, \ \ \ 
P=P_j\circ W_j, \ \ \
Q_j=V_j\circ W_j, \ \ \  \ \ \ W_j(a)=W_j(b).
$$ Moreover,
one of the following conditions holds:
\vskip 0.25cm 
\noindent 1) $r=2$ and 
$$P=U\circ x^{2}R^2(x^2) \circ V,\ \ \ W_1=x^2 \circ V, \ \ \ W_2=xR(x^2)\circ V,$$ 
where $R,$ $U$, $V$  are real polynomials;
\vskip 0.25cm
\noindent 2)  $r=2$ and
$$P=U\circ T_{m_1m_2} \circ V,\ \ \ W_1= T_{m_1} \circ V, \ \ \  
W_2= T_{m_2}\circ V,$$ 
where $U$, $V$ are  real polynomials, $m_1>1,$ $m_2>1$, $\GCD(m_1,m_2)=1;$
\vskip 0.25cm 
\noindent 3) $r=3$ and 
$$P=U\circ x^{2}R^2(x^{2})\circ  T_{m_1m_2}
\circ V,$$ 
$$W_1=T_{2m_1} \circ V, \ \ \ W_2=T_{2m_2}\circ V, \ \ \ 
W_3= (xR(x^2) \circ T_{m_1m_2}) \circ V,$$ 
where $R,$ $U$, $V$ are  real polynomials, 
$m_1>1,$ $m_2>1$ are odd, and $\GCD(m_1,m_2)=1.$
\et
\pr Our strategy is to apply Theorem \ref{mpc} and to use the condition that $P,Q\in \R[x]$ and $a,b\in \R.$ Assume first that \eqref{c} holds for some $\t P, \t Q, W\in \C[z].$
Applying Corollary \ref{reco} to the equality $P=\t P\circ W$ we conclude  that without loss of generality we may assume that $\t P$ and $W$ are contained in $\R[x]$.  Now the equality $Q=\t Q\circ W$ 
implies by  Corollary \ref{reco0} that 
$\t Q$ also is contained in $\R[x]$.

Assume that $Q$ is a solution of the first type. It follows from the equality $P=P_1\circ W_1$ by Corollary \ref{reco} that there exists a complex polynomial of degree one $\mu_1$ such that the polynomial $\mu_1\circ W_1$ has real coefficients. Further, applying
Corollary \ref{reco}  to 
the equality $\mu_1\circ W_1=\mu_1\circ z^n\circ V$, we conclude that  there exists a complex polynomial of degree one $\mu_2$ such that the polynomials $\mu_1\circ z^n\circ \mu_2$ and $\mu_2^{-1}\circ V$ have real coefficients. By Lemma \ref{eche}, a) this implies that 
there exist $\t \mu\in \R[x]$ and $c\in \C$ such that 
$\mu_2=c\t \mu.$ Since $\mu_2^{-1}=\t\mu^{-1}\circ z/c,$ it follows now from $\mu_2^{-1}\circ V\in \R[x]$ that $V/c\in \R[x].$ Therefore, changing the polynomial $V$ to $V/c$, and modifying the polynomials 
$P_1,$ $V_1,$ and $R$ in an obvious way, without loss of generality we may assume that $V\in \R[x].$ 

Clearly, $V\in \R[x]$ implies that $W_1=z^n\circ V\in \R[x]$. It follows now from  $P=P_1\circ W_1$
by Corollary \ref{reco0} that $P_1\in \R[x]$. Furthermore, 
it follows from $W_1(a)=W_1(b)$ and $a,b\in \R$ that $n=2k$, and $V(a)=-V(b).$
Since $\GCD(s,n)=1,$ the polynomial $z^sR(z^n)$ has the form $z\t R(z^2)$ for some $\t R\in \C[z]$.  
Thus, changing $P_1$ tp $P_1\circ z^{n/2}$ and $z^sR(z^n)$ to $z\t R(z^2)$, 
without loss of generality we may assume that $W_1=z^2\circ V$ and 
$W_2=zR(z^2)\circ V$. Applying 
Corollary \ref{reco0} to the equality $P=(P_2\circ zR(z^2))\circ V$ we see that $P_2\circ zR(z^2)\in \R[x]$. Therefore, taking into account that the constant term of $zR(z^2)$ is zero, 
Corollary \ref{reco0} implies that for $c\in \C$ such that the leading coefficient of 
$czR(z^2)$ is real the polynomials $P_2\circ z/c$ and  $czR(z^2)$ are contained in $\R[x]$. 
Thus, modifying the polynomials $P_2$ and $R$ we can assume that they are  
contained in $\R[x]$. Now Corollary \ref{reco0} applied to the equality $P=U\circ (z^2R^2(z^2)\circ V)$ implies that $U\in \R[x]$.

Finally, the equality 
$$Q=V_1\circ W_1+V_2\circ W_2=\overline{V}_1\circ W_1+\overline{V}_2\circ W_2$$ implies that 
$$Q=\frac{V_1+\overline{V}_1}{2}\circ W_1+\frac{V_2+\overline{V}_2}{2}\circ W_2.$$
Therefore, changing if necessary $V_1$ to $(V_1+\overline{V}_1)/2$ and $V_2$ to $(V_2+\overline{V}_2)/2,$ without loss of generality we may assume that $V_1,$ $V_2\in \R[x].$

Assume now that $Q$ is a solution of the third type. We may assume that $V(a)\neq - V(b)$, for otherwise, as it was observed after Theorem \ref{mpc}, this solution also belongs to the first type considered earlier.
As above, we conclude that there
exist complex polynomials of degree one $\mu_1$ and $\mu_2$ such that the polynomials $\mu_1\circ T_{2m_1}\circ \mu_2$ and $\mu_2^{-1}\circ V$ have real coefficients. By Lemma \ref{eche}, b), this implies that 
either $\mu_2\in \R[x]$ or 
there exists $\t \mu\in \R[x]$ such that 
$\mu_2=i\t \mu.$ Since, $\mu_2^{-1}\circ V\in \R[x],$ in the first case $V\in \R[x],$ while in the second one, \be \l{asdf} V=i\t V, \ \ \ \t V\in \R[x].\ee Let us show that equality \eqref{asdf} is impossible. 
Indeed, applying Lemma \ref{skun}, a) to the equalities $W_1(a)=W_1(b),$ $W_2(a)=W_2(b)$ and arguing as in Lemma \ref{red}, we conclude that the numbers $V(a)$ and $V(b)$ are roots of the polynomial $T_{4m_1m_2}^{\prime}$, for otherwise $V(a)= - V(b)$. Since 
$T^{\prime}_n$ has only real zeroes, we conclude that $V(a),V(b) \in \R$, and hence \eqref{asdf} is impossible in view of $a,b\in \R.$
Thus, $V\in \R[x]$.

Applying now Corollary \ref{reco0} to the equality $P=(P_3\circ zR(z^2))\circ (T_{m_1m_2}\circ V)$ we conclude that $P_3\circ zR(z^2)\in \R[x]$. Furthermore, arguing as above, we conclude that without loss of generality we may assume that 
$P_3,R\in \R[x]$ as well as $P_1,P_2,U\in \R[x]$ and $V_1,V_2,V_3\in \R[x].$

The 
proof of the theorem in the case where $Q$ is a solution of the second type is obtained similarly with obvious simplifications. \qed

\section{Proof of Theorem \ref{t2}}
\subsection{Plan of the proof}
In the rest of the paper we always will assume that all considered polynomials  have real coefficients.

Let us describe a general plan of the proof of Theorem \ref{t2}. First, observe that we may assume that 
\be \la{lu} \R(P,Q)=\R(x).\ee Indeed, otherwise by the L\"uroth theorem, $\R(P,Q)=\R(W)$ for some $W\in \R(x),$ $\deg W\geq 2,$ implying that 
$$P=\t P \circ W, \ \ \ Q=\t Q\circ W$$ for some $\t P, \t Q \in \R(x)$ such that $\R(\t P,\t Q)=\R(x).$ Moreover, since $P,Q\in \R[x]$, it is easy to see that 
we may assume that $\t P, \t Q, W\in\R[x].$ Therefore, since equalities  \eqref{mix} imply the equalities  
\be  
\int_{W(a)}^{W(b)} \t P^id \t Q=0, \ \ \  \int_{W(a)}^{W(b)} \t Q^jd \t P=0, \ \ \ i,j \geq 0,
\ee 
it is enough to prove that if  \eqref{lu} holds, then $P,Q$ cannot satisfy \eqref{mix} for  $a\neq b.$ 

Applying Theorem \ref{mpr} to the first and to the second system of equations in \eqref{mix} separately, we arrive to nine different
``cases'' depending on types of solutions appearing in Theorem \ref{mpr}. For example, ``the case (2,1)'' means that $Q$ is a solution of the second type 
of the polynomial moment problem \eqref{mom}, while $P$ is a solution of the first type of the polynomial moment problem
$$\int_a^b Q^jd P=0, \ \ \ j \geq 0.$$
In more details, this means that, from one hand, 
$$Q=V_1\circ  W_1+ V_2\circ  W_2,\ \ \ P=U\circ T_{nm} \circ  V,$$ where 
\be \la{uio}  W_1=T_n \circ  V, \ \ \ W_2=T_m\circ  V,\ee and 
$$W_1(a)=W_1(b), \ \ \ W_2(a)=W_2(b),$$ while, from the other hand,
$$P=\t V_1\circ \t W_1+\t V_2\circ \t W_2 ,\ \ \ Q=\t U\circ  x^{2} R^2(x^2)\circ\t V,$$ where 

$$\t W_1= x^2 \circ \t V, \ \ \  
\t W_2=xR(x^2)\circ \t V,$$  and 
$$\t W_1(a)=\t W_1(b), \ \ \ \t W_2(a)=\t W_2(b).$$
In view of assumption \eqref{lu}, the polynomial $V$ (as well as the polynomial $\t V$) is of degree one
for otherwise
$$\R(P,Q)\subseteq\R(V)\subsetneq \R(x).$$ Furthermore, it is clear that without loss of generality we may assume that one of the polynomials $V$ and $\t V$ equals $x.$
Our strategy will be to show that such systems of equations always imply that equalities \eqref{c} hold, in contradiction with \eqref{lu} (recall that 
the condition $W(a)=W(b)$  implies that $\deg W>1$).

Since we may exchange $P$ and $Q$, it is necessary to consider only the cases (1,1), (2,1), (2,2), (3,1), (3,2), and (3,3).   
Finally, 
we may impose some additional restrictions related to the fact that  a solution  of the (usual) polynomial moment problem may belong to different types. For example, assuming that the theorem is already proved in the case (1,1), considering the case (2,1) we may assume that $n>2,$ $m>2$ in \eqref{uio}, since otherwise the solution 
$P,Q$  also  belongs to the case (1,1).

For a polynomial $$P=a_nx^n+a_{n-1}x^{n-1}+ ... +a_1x+a_0, \ \ \ a_i\in \R, \ \  \ 0\leq i \leq n,$$ of degree $n$, set $$C_i(P)=a_{n-i}, \ \ \ 0\leq i \leq n.$$

The following simple lemma permits to control initial terms in a composition 
of two polynomials and is widely used in the following. 

\bl \la{l1} Let $T$ be a polynomial of degree $d$.
Then for any polynomial $S$ of degree $r$ with the leading coefficient $c$ the equalities  
\be \la{wer} C_i(S\circ T)=C_i(cx^r\circ T), \ \ \ 0\leq i \leq d-1,\ee hold. In particular, 
for any two polynomials $S_1,$ $S_2$ of equal degree  with equal leading coefficients
the equalities  
$$ C_i(S_1\circ T)=C_i(S_2\circ T), \ \ \ 0\leq i \leq d-1,$$ hold.   
\el  
\pr Indeed, $\deg (S-cx^r)\circ T=dr-d$. Therefore, \eqref{wer} holds. \qed

\bc \la{c1} Let $T(z)$ be a polynomial of degree $d\geq 2$  such that  $C_1(T)=0$ holds,
$U$ be an arbitrary polynomial, and $\alpha, \beta\in \R,$ $\alpha\neq 0.$ Then  equality $C_1(U\circ T\circ (\alpha x+\beta))=0$ holds if and only if $\beta=0.$ 
\ec
\pr Indeed,   if $\deg U=r$, $C_0(U)=c$, then 
$C_1(U\circ T)=C_1(cz^r\circ T)$ by Lemma \ref{l1}. On the other hand, $C_1(cz^r\circ T)=0$, since $C_1(T)=0$. Therefore, $C_1(U\circ T)=0$, and it is easy to see that for any polynomial $F$ such that  $C_1(F)=0,$ the equality 
$C_1(F\circ (\alpha x+\beta))= 0$  holds if and only if $\beta=0.$ 
\qed

\subsection{Proof of Theorem \ref{t2} in the cases (1,1)}

\bl \la{ll2} Let $W_1$, $W_2$ be polynomials of degree two such that 
$W_1(a)=W_1(b),$ $W_2(a)=W_2(b)$ for distinct $a,b\in \R$. Then $W_2=\lambda_1W_1+\lambda_2$ for some $\lambda_1,$ $\lambda_2 \in \R.$ 
\el
\pr 
Let 
$$P=\alpha_1x^2+\beta_1x+\gamma_1, \ \ \ \ Q=\alpha_2x^2+\beta_2x+\gamma_2,$$ 
where $\alpha_1,\beta_1,\gamma_1,\alpha_2,\beta_2,\gamma_2\in \R.$
Then conditions of the lemma yield the equalities 
$$\alpha_1(a+b)+\beta_1=0, \ \ \ \alpha_2(a+b)+\beta_2=0,$$ implying the statement.
\qed
In order to proof the theorem in the case (1,1) it is enough to observe that in this case there exist $U,R,\t U, \t R\in \R[x]$ such that  
$$P=(U\circ zR(z^2))\circ W_1, \ \ \ Q=(\t U\circ z\t R(z^2))\circ \t W_1$$ 
where 
$$W_1=z^2\circ V, \ \ \ \ \t W_1= z^2\circ \t V,$$ 
are polynomials of degree two such that 
$W_1(a)=W_1(b),$ $\t W_1(a)=\t W_2(b)$. By Lemma \ref{ll2}, we have
$\t W_2=\lambda_1W_1+\lambda_2$ for some $\lambda_1,$ $\lambda_2 \in \R$, and hence  
\eqref{c} holds for $W=W_1.$

\subsection{Proof of Theorem \ref{t2} in the case (2,1)}
If $P,Q$ is a solution of \eqref{mix} corresponding to the case (2,1), then without loss of generality we may assume that 
there exist polynomials $V_1,V_2, U, \t V_1,\t V_2, \t U,  R$ and $\alpha,\beta\in \R,$ $\alpha\neq 0,$ such that 
\be \la{urab} Q=V_1\circ T_{m_1}+V_2\circ T_{m_2}
=\t U\circ x^2{R}^2(x^2)\circ (\alpha x+\beta),\ee and 
$$ P=\t V_1\circ x^2 \circ (\alpha x+\beta)+\t V_2\circ xR(x^2) \circ (\alpha x+\beta)
=U\circ T_{m_1m_2} ,$$
where  $\GCD(m_1,m_2)=1.$ In addition, for 
the polynomials $$ W_1=T_{m_1}, \ \ W_2=T_{m_2}, \ \ \t W_1=z^2\circ (\alpha x+\beta),\ \  \t W_2=xR(x^2)\circ (\alpha x+\beta)$$
the equalities \be \la{ko2} W_1(a)=W_1(b), \ \ \ W_2(a)=W_2(b), \ \ \ \t W_1(a)=\t W_1(b), \ \ \ \t W_2(a)=\t W_2(b)\ee hold.
If at least one of the numbers $m_1,m_2$ equals 2, then $P,Q$ belongs to the type (1,1) considered above.
So, we may assume that  $m_1\geq 3,$ $m_2\geq 3$. Notice that the second representation for $Q$ in \eqref{urab} implies that $n=\deg Q$ is even and the inequality $n\geq 6$ holds.

Although the above conditions seem to be very strong, it is difficult to use them in their full generality since they contain many unknown parameters. Thus, actually we mostly will use only the fact that the right part of  \eqref{urab} is a polynomial in $x^2\circ (\alpha x+\beta)$ together with first three equalities in \eqref{ko2}.

First of all observe that any polynomial of the form $Q=V_1\circ T_{m_1}+V_2\circ T_{m_2}$ can be represented in the form  \be \la{xy} Q=d_nT_n+d_{n-1}T_{n-1}+\dots +d_{1}T_1+d_0,\ \ \ \ \ d_i\in \R,\ee
where $d_i=0$, unless $i$ is divisible either by $m_1$ or by $m_2$.
Indeed, it is clear that $T_0,T_1, \dots , T_{r}$ is a basis of a subspace of $\R[x]$ consisting of all polynomials of degree $\leq r$. 
Therefore, a polynomial $P$ can be represented in the form 
$$P=V\circ T_m=(a_rx^r+a_{r-1}x^{r-1}+...+a_1x+a_0)\circ T_{m}$$ if and only if 
$$P=(b_rT_r+b_{r-1}T_{r-1}+ ... +b_1T_1+b_0)\circ T_m=b_rT_{rm}+b_{r-1}T_{(r-1)m}+ ... +b_1T_m+b_0.$$ Consequently,  
$Q=V_1\circ T_{m_1}+V_2\circ T_{m_2}$ can be represented in the required form. 

Define $C(n,m_1,m_2)$ as the set of all polynomials \eqref{xy} such that $d_i=0$, unless $i$ is divisible either by $m_1$ or by $m_2$, and $d_n\neq 0.$
To be definite, we always will assume that $n$ is divisible by $m_1.$
Similarly to the notation $C_i(Q)$ introduced above, for a polynomial $Q\in C(n,m_1,m_2)$ set 
$$Ch_i(Q)=d_{n-i}, \ \ \ 0\leq i \leq n.$$

\bl \la{dfg} Let $Q\in C(n,m_1,m_2)$, where $m_1$ and $m_2$ are 
coprime, $m_1\geq 3,$ $m_2\geq 3$. 

\vskip 0.2cm
\noindent 1) If $Ch_1(Q)\neq 0$, then $Ch_2(Q)=0.$
\vskip 0.2cm
\noindent 2) If $Ch_1(Q)\neq 0,$ $Ch_3(Q)\neq 0,$ then $m_1=3$, and $Ch_4(Q)= 0,$  
\vskip 0.2cm
\noindent 3)  If $Ch_1(Q)\neq 0,$ $Ch_3(Q)\neq 0,$ $Ch_5(Q)\neq 0,$ then $m_1=3,$ $m_2=4.$
\el 
\pr 
If $Ch_1(Q)\neq 0$, then  $m_2\mid n-1.$  
Further, clearly 
\be \la{inv} \vert m_ik_1-m_ik_2\vert\geq m_i\geq 3, \ \ \ i=1,2,\ \ \ k_1,k_2\in \N,\ee   unless $k_1=k_2$,
implying that 
$n-2$ may not be divisible neither by $m_1$, since $m_1\mid n$, 
nor  by $m_2$, since $m_2\mid n-1$. 
Thus, $Ch_2(Q)=0.$

Assume that additionally $Ch_3(Q)\neq 0$. Since $m_2\mid n-1$, the number $n-3$ may not be divisible by $m_2$ in view of 
\eqref{inv}. Therefore, $n-3$ is divisible by $m_1$ implying that $m_1=3$. Furthermore, $Ch_4(Q)=0$ for otherwise \eqref{inv} implies that $m_2=3$ in contradiction with $\GCD(m_1,m_2)=1.$ 

Finally, if additionally $Ch_5(Q)\neq 0$, it follows from $m_1=3$ by \eqref{inv}
that $m_2\vert n-5$ implying that $m_2=4.$ \qed

\bc \la{ddfgg} Let $Q\in C(n,m_1,m_2)$, where $m_1$ and $m_2$ are 
coprime, $m_1\geq 3,$ $m_2\geq 3$, and $n$ is even. If $Ch_1(Q)\neq 0$, then 
either $Ch_2(Q), Ch_3(Q)$ vansih, or $Ch_2(Q)$, $Ch_4(Q),$ $Ch_5(Q)$ vanish.

\ec
\pr Since $m_2\vert n-1$, and $n$ is even, $m_2$ is odd. Therefore, $m_2\neq 4$, and the statement follows from  
Lemma \ref{dfg}. \qed

\bl \la{l2} Let $Q=F\circ x^2\circ (x-\delta)$, $\delta \in \R,$ where $n=\deg Q\geq 6.$ 
\vskip 0.2cm
\noindent 1) If the coefficients $Ch_2(Q)$, $Ch_3(Q)$ vanish, then either $\delta=0,$ or \be \la{azx} (2\delta)^2=\frac{3}{(n-1)(n-2)}. \ee
\vskip 0.2cm
\noindent 2) If the coefficients $Ch_2(Q)$, $Ch_4(Q),$ $Ch_5(Q)$ vanish, then either $\delta=0,$ or $2\delta$ satisfies the equation 
\be \la{xza} \frac{2}{15}(n-1)(n-2)(n-3)(n-4)t^4-(n-2)(n-3)t^2+1=0.\ee
\el 
\pr If $Ch_2(Q)$, $Ch_3(Q)$ vanish, then $Q=c_0T_n+c_1T_{n-1}+ R_{1},$ where $\deg R_1\leq n-4$.
Set $$T_s^*(x)=2T_s\left(\frac{x}{2}\right), \ \ \ s\geq 1.$$ Clearly, the equality 
$$ c_0T_n+c_1T_{n-1}+ R_{1}=F\circ x^2\circ (x-\delta)$$ implies the equality 
\be \la{0100} c_0T_n^*+c_1T_{n-1}^*+ \widetilde R_{1}=\t F\circ x^2\circ (x-\gamma),\ee
where $\t R_1=2R_1(x/2),$ $\t F=2F(x/4),$ and
$\gamma=2\delta.$
Furthermore, without loss of generality we may assume that $c_0=1.$

Since the right part of \eqref{0100} is a polynomial in $(x-\gamma)^2$, it follows from  \eqref{0100} taking into account the inequality $\deg R_1\leq n-4$
that the derivatives of $T_n^*+c_1T_{n-1}^*$ of orders $n-1$ and $n-3$ at the point $\gamma$ vanish, that is  
$$T_n^{*(n-1)}(\gamma)+c_1T_{n-1}^{*(n-1)}(\gamma)=0,$$
\vskip -0.2cm
$$T_n^{*(n-3)}(\gamma)+c_1T_{n-1}^{*(n-3)}(\gamma)=0.$$ 
Since 
\be \la{f} T_s^*=x^s-sx^{s-2}+\frac{s(s-3)}{2}x^{s-4}-\frac{s(s-4)(s-5)}{6}x^{s-6}+ ... \ee by \eqref{cheb}, this implies that
$$n!\gamma+c_1(n-1)!=0,$$
$$\frac{n!}{3!}\gamma^3-n(n-2)!\gamma+c_1\left(\frac{(n-1)!}{2!}\gamma^2-(n-1)(n-3)!\right)=0.$$
The first of these equalities implies 
that $c_1=-n\gamma$. Substituting this value of $c_1$ in the second equality we obtain
$$\frac{n!}{3!}\gamma^3-n(n-2)!\gamma-n\gamma\left(\frac{(n-1)!}{2!}\gamma^2-(n-1)(n-3)!\right)=-\frac{n!}{3}\gamma^3+n(n-3)!\gamma=0$$ implying that $2\delta= \gamma$ satisfies \eqref{azx}, unless $\delta=0.$

Similarly, if $Ch_2(Q)$, $Ch_4(Q),$ $Ch_5(Q)$ vanish, we arrive to the equality 
\be \la{1} T_n^*+c_1T_{n-1}^*+c_3T_{n-3}^*+ \widetilde R_{1}=\t R\circ x^2\circ (x-\gamma),\ee
where $\deg \t R_1\leq n-6$ and $\gamma=2\delta,$ implying 
that the derivatives of $T_n^*+c_1T_{n-1}^*+c_3T_{n-3}^*$ of orders $n-1$ and $n-3,$ and $n-5$ at the point $\gamma$ vanish. Thus, 
$$T_n^{*(n-1)}(\gamma)+c_1T_{n-1}^{*(n-1)}(\gamma)=0,$$
\vskip -0.2cm
$$T_n^{*(n-3)}(\gamma)+c_1T_{n-1}^{*(n-3)}(\gamma)+c_3T_{n-3}^{*(n-3)}(\gamma)=0,$$

$$T_n^{*(n-5)}(\gamma)+c_1T_{n-1}^{*(n-5)}(\gamma)+c_3T_{n-3}^{*(n-5)}(\gamma)=0,$$
or equivalently 
$$n!\gamma+c_1(n-1)!=0,$$
$$\frac{n!}{3!}\gamma^3-n(n-2)!\gamma+c_1\left(\frac{(n-1)!}{2!}\gamma^2-(n-1)(n-3)!\right)+c_3(n-3)!=0,$$ and
$$\frac{n!}{5!}\gamma^5-\frac{n(n-2)!}{3!}\gamma^3+\frac{n(n-3)(n-4)!}{2}\gamma+$$ $$+c_1\left(\frac{(n-1)!}{4!}\gamma^4 -\frac{(n-1)(n-3)!}{2!}\gamma^2+\frac{(n-1)(n-4)(n-5)!}{2}\right)+$$
$$+c_3\left(\frac{(n-3)!}{2!}\gamma^2-(n-3)(n-5)!\right)=0.$$
As above, it follows from the first of these equalities by \eqref{f} that $c_1=-n\gamma$ and substituting this value of $c_1$ in the second equality we obtain
$$\frac{n!}{3!}\gamma^3-n(n-2)!\gamma-n\gamma\left(\frac{(n-1)!}{2!}\gamma^2-(n-1)(n-3)!\right)+c_3(n-3)!=0$$ implying that 
$$c_3=\frac{n(n-1)(n-2)}{3}\gamma^3-n\gamma.$$
Now the third equality gives us 
$$\frac{n!}{5!}\gamma^5-\frac{n(n-2)!}{3!}\gamma^3+\frac{n(n-3)(n-4)!}{2}\gamma-n\gamma\left(\frac{(n-1)!}{4!}\gamma^4 -\frac{(n-1)(n-3)!}{2!}\gamma^2+\frac{(n-1)(n-4)(n-5)!}{2}\right)+$$
$$+\left(\frac{n(n-1)(n-2)}{3}\gamma^3-n\gamma\right)\left(\frac{(n-3)!}{2!}\gamma^2-(n-3)(n-5)!\right)=0.$$
The coefficient of $\gamma^5$ in this expression is 
$$\frac{n!}{3!}\left(\frac{1}{20}-\frac{1}{4}+1\right)=\frac{2n!}{15}\ .$$
The coefficient of $\gamma^3$  is 
$$-\frac{n(n-2)!}{3!}+\frac{n(n-1)(n-3)!}{2!}-\frac{n(n-1)(n-2)(n-3)(n-5)!}{3}-\frac{n(n-3)!}{2!}=$$
$$= -\frac{n(n-3)!}{2!}\left(\frac{n-2}{3}-(n-1)\right)-n(n-3)(n-5)!\left(\frac{(n-1)(n-2)}{3}+\frac{(n-4)}{2}\right)=$$
$$=\frac{(2n-1)n(n-3)!}{6}-n(n-3)(n-5)!\left(\frac{2n^2-3n-8}{6}\right)=$$ $$=\frac{n(n-3)(n-5)!}{6}\left((2n-1)(n-4)-(2n^2-3n-8)\right)=-n(n-2)(n-3)(n-5)!$$

Finally, the coefficient of $\gamma$ is 
$$\frac{n(n-3)(n-4)!}{2}-\frac{n(n-1)(n-4)(n-5)!}{2}+n(n-3)(n-5)!=
$$
$$=n(n-5)!\left(\frac{(n-3)(n-4)}{2}-\frac{(n-1)(n-4)}{2}+(n-3)\right)=n(n-5)!\, . 
$$
Collecting terms and canceling by $n(n-5)!$ we see that $2\delta=\gamma$ satisfies  
\eqref{xza}, unless $\delta=0.$ \qed

\bc \la{lal} Let $Q=R\circ x^2\circ (x-\delta)$, $\delta \in \R.$
\vskip 0.2cm
\noindent 1) If the coefficients $Ch_2(Q)$, $Ch_3(Q)$ vanish and $n=\deg Q\geq 6,$ then the number $4 \delta$ 
is not an  algebraic integer, unless $\delta=0.$

\vskip 0.2cm
\noindent 2) If the coefficients $Ch_2(Q)$, $Ch_4(Q),$ $Ch_5(Q)$ vanish and $n=\deg Q\geq 9,$ 
then the number $4 \delta$ 
is not an  algebraic integer, unless $\delta=0.$
\ec 
\pr 
Set $\gamma=4\delta.$ If $Ch_2(Q)$, $Ch_3(Q)$ vanish and $\delta\neq0,$ then $\gamma$ is a root of the equation 
$$t^2-\frac{12}{(n-1)(n-2)}=0.$$ 
Since for $n\geq 6$ the number $\frac{12}{(n-1)(n-2)}$ is not an integer, this implies that  $\gamma$ cannot be an  algebraic integer
of degree two. Moreover, $\gamma$ cannot be an  algebraic integer
of degree one for otherwise $\gamma$ is an integer implying that $\gamma^2=\frac{12}{(n-1)(n-2)}$ is also an integer.

If $Ch_2(Q)$, $Ch_4(Q),$ $Ch_5(Q)$ vanish and $\delta\neq0,$ then $\gamma$ is a root of the equation  
\be \label{ur} (n-1)(n-2)(n-3)(n-4)t^4-30(n-2)(n-3)t^2+120=0.\ee
If this equation is irreducible over $\Q$, then $\gamma$ cannot be an algebraic integer since 
the number $(n-1)(n-2)(n-3)(n-4)$ obviously does not divide 120 for $n\geq 9.$  

Assume now that $\gamma$ is an algebraic integer satisfying an irreducible equation $t^2+c_1t+c_2=0,$ $c_1,c_2\in \Z.$ Then 
by the Gauss lemma the equality 
$$(n-1)(n-2)(n-3)(n-4)t^4-30(n-2)(n-3)t^2+120=(t^2+c_1t+c_2)(d_0t^2+d_1t+d_2)$$ holds for some $d_0,d_1,d_2\in \Z$.
Since the coefficients of  $t^3$ and $t$ in the left part vanish we have:
$$d_1+c_1d_0=0, \ \ \ \ c_1d_2+c_2d_1=0,$$ implying that \be \la{bege1} d_1=-c_1d_0, \ \ \ d_2=c_2d_0,\ee unless
\be \la{bege2} c_1=0,\ \ \  d_1=0.\ee

If \eqref{bege1} holds,  then 
$120=c_2d_2=c_2^2d_0$ in contradiction with 
$d_0=(n-1)(n-2)(n-3)(n-4)$ and $n\geq 9.$
Similarly, if \eqref{bege2} holds, 
then $\gamma^2=-c_2$ is an integer, and  
it follows from  \eqref{ur} that the number
$(n-2)(n-3)\gamma^2$ divides 120 implying easily a contradiction with the condition $n\geq 9.$ 

Finally, observe that if $\gamma$ is a rational root of $\eqref{ur}$, then $-\gamma$ also is a root of \eqref{ur} and $t^2-\gamma^2$ is a divisor of 
\eqref{ur}. Therefore, any irreducible over $\Q$ factor of \eqref{ur} has the degree one, two, or four. Thus, in order to finish the proof it is enough 
to observe 
that $\gamma$ is not an integer, since otherwise $\gamma^2$ is also an integer implying as above that 
$(n-2)(n-3)\gamma^2$ divides 120. \qed

\noindent{\it Proof of Theorem \ref{t2} in the case (2,1).} 
Observe first that if $\beta=0$ in \eqref{urab}, then the theorem is true. Indeed, in this case the condition $\t W_1(a)=\t W_1(b)$ is equivalent 
to the condition $T_2(a)=T_2(b).$
Therefore, applying Lemma \ref{skun}, b) to the equalities 
$$T_{m_1}(a)=T_{m_1}(b), \ \ \ T_{m_2}(a)=T_{m_2}(b)\ \ \ T_2(a)=T_2(b),$$
we conclude that at least one of the numbers $m_1,$ $m_2$ is even and hence $Q=\t U\circ xR^2(x)\circ (\alpha x)^2 $ and 
$P=U\circ T_{m_1m_2}=U\circ T_{m_1m_2/2}\circ T_2$ satisfy \eqref{c} for $W=x^2.$

Further, observe that \eqref{urab} implies that 
\be\la{svina} Q=\t U\circ x{R}^2(x)\circ x^2\circ (\alpha x+\beta)=F\circ x^2 \circ \left(x+\frac{\beta}{\alpha}\right),\ee
where $F=\t U\circ x{R}^2(x)\circ \alpha^2 x,$ while the condition $\t W_1(a)=\t W_1(b)$ yields that 
\be \la{a+b} a+b=-\frac{2\beta}{\alpha}.\ee

If $Ch_1(Q)=0$, then \eqref{xy} implies that also $C_1(Q)=0$, since $C_1(T_n)=0$ by \eqref{cheb}. 
In its turn $C_1(Q)=0$ implies that $\beta=0$ by Corollary \ref{c1} applied to \eqref{svina}. So, 
assume that $Ch_1(Q)\neq 0$. By Corollary \ref{ddfgg}, this implies that
either $Ch_2(Q), Ch_3(Q)$ vanish, or $Ch_2(Q)$, $Ch_4(Q),$ $Ch_5(Q)$ vanish.

If $Ch_2(Q), Ch_3(Q)$ vanish, then 
Corollary \ref{lal}, 1)  applied to \eqref{svina}  implies that  the number $-4 \beta/\alpha$ 
if  not an  algebraic integer, unless $\beta=0$. On the other hand, \eqref{a+b} implies that 
$-4 \beta/\alpha$ is an  algebraic integer, since $2a$ and $2b$ are algebraic integers by Corollary \ref{red}.
Thus, we conclude again that $\beta=0.$

Similarly, assuming that $Ch_3(Q)\neq 0$, while
$Ch_2(Q)$, $Ch_4(Q),$ $Ch_5(Q)$ vanish, 
we may apply Corollary \ref{lal}, 2), 
whenever $n\geq 9.$ Therefore, 
since $Ch_3(Q)\neq 0$ implies that $3\vert n$ in view of the equality $m_1=3$, 
the only case 
which remains uncovered is the one where $n=6$. In this case the inequality $Ch_1(Q)\neq 0$ implies that $m_2=5.$  
Notice that it follows from \eqref{ur} that for $n=6$ the number  $4 \beta/\alpha$ satisfies the equation $t^4-3t^2+1=0$ whose roots
are algebraic integers. Moreover, the system 
$$T_5(a)=T_5(b), \ \ T_3(a)=T_3(b), \ \  (2a+2b)^4-3(2a+2b)^2+1=0$$ 
has non-zero solutions. Thus, for $n=6$ the previous reasoning fails.

In order to prove the theorem in this case remind that in Lemma \ref{l2} we used a condition which is weaker than the one in \eqref{urab}. Therefore, in order to finish the proof 
it is enough to show that 
the equality  
\be \la{tgv} T_6+c_1T_5+c_3T_3+c_6=c(z(z^2-d))^2\circ (z-\beta),\ee
where $c_1,c_3,c_6,c,d,\beta\in \R,$ is possible only if $c_1=0$. 
This statement may be verified by 
a direct calculation. Namely, the comparison of leading coefficients of both parts of \eqref{tgv} implies that $c=32,$ while the comparison of other coefficients gives 
$$16 c_1+192 b=0, \ \ \  -480 b^2+64 d-48=0,$$ $$640 b^3-256 b d-20 c_1+4c_3=0, \ \ \ -480 b^4+384 b^2 d-32 d^2+18=0,$$ $$192 b^5-256 b^3 d+64 b d^2+5 c_1-3 c_3,  \ \ \  c_6=-32 b^6+64 b^4 d-32 b^2 d^2-1.
$$
We leave the reader to verify (for example, with the help of Maple) that the only solution of the above system is 
$$
c_1 = 0,\ \  c_3 = 0,\  \ c_6 = 1, \  \ b = 0, \ \ d = 3/4,$$
(for these values of parameters equality 
\eqref{tgv} simply reduces to the equality $T_6=T_2\circ T_3$).

\subsection{Proof of Theorem \ref{t2} in the case (2,2)}
First, observe that Theorem \ref{t2} in the case (2,2) follows from the following statement.

\bp \la{vot} Let $V_1,V_2, U, \t V_1, \t V_2, \t U\in \R[x]$ and $\alpha,\beta\in \R$, $\beta\neq 0,$ satisfy the equalities 
\be \la{uraba} Q=V_1\circ T_{m_1}+V_2\circ T_{m_2}=\t U\circ T_{\t m_1\t m_2}\circ (\alpha x+\beta),\ee and 
\be \la{uraba2} P=\t V_1\circ T_{\t m_1}\circ (\alpha x+\beta)+\t V_2\circ T_{\t m_2}\circ (\alpha x+\beta)=U\circ T_{m_1m_2},\ee 
where \be \la{kvak} \GCD(m_1,m_2)=1, \ \ \ \GCD(\t m_1,\t m_2)=1,\ee and $m_1\geq 3, m_2\geq 3, \t m_1\geq 3,\t m_2\geq 3.$ 
Then $\alpha=\pm 1,$ $\beta =0.$ 

\ep 

Indeed, in the case (2,2) conditions \eqref{uraba}, \eqref{uraba2}, \eqref{kvak} are satisfied for some $m_1\geq 2,$ $m_2\geq 2,$ $ \t m_1\geq 2,$ $\t m_2\geq 2.$ 
Additionally, the equalities 
\be \la{asdd-} T_{m_1}(a)=T_{m_1}(b), \ \ \ T_{m_2}(a)=T_{m_2}(b),\ee 
\be \la{asdd} T_{\t m_1}(\alpha a+\beta)=T_{\t m_1}(\alpha b+\beta), \ \ \ 
T_{\t m_2}(\alpha a+\beta)=T_{\t m_2}(\alpha b+\beta)\ee
hold.
If at least one of the numbers $m_1,m_2$ equals 2, then $P,Q$ belongs to the type (1,2) considered above. So, we may assume that  $m_1\geq 3,$ $m_2\geq 3$. Similarly, we may assume that   
$\t m_1\geq 3,\t m_2\geq 3,$ since otherwise $P,Q$ belongs to the type (2,1).
Since under these conditions Proposition \ref{vot} implies that 
$\alpha=\pm 1,$ $\beta =0,$ it follows from the equalities 
$$T_{m_1}(a)=T_{m_1}(b), \ \ \ T_{m_2}(a)=T_{m_2}(b), \ \ \ T_{{\t m}_1}(\pm a)=T_{{\t m}_1}(\pm b)$$ 
by  Lemma \ref{skun}, b), taking into account \eqref{kots}, that 
$T_s(a)=T_s(b)$, where either $s=\GCD(m_1, {\t m}_1)$, or $s=\GCD(m_2, {\t m}_2).$
Since
$$Q=\t U\circ T_{\t m_1\t m_2}\circ (\pm z)=\t U\circ (\pm  T_{\t m_1\t m_2})=\t U\circ (\pm  T_{\t m_1\t m_2/s})\circ T_s$$
and 
$$P=U\circ  T_{m_1m_2}=U\circ  T_{m_1m_2/s}\circ T_s,$$
we conclude that \eqref{c} holds for $W=T_s.$

\vskip 0.2cm

The next lemmas are similar to Lemma \ref{dfg} and Lemma \ref{l2} and are used for imposing restrictions on possible values of $\alpha$ and $\beta$ in \eqref{uraba}, \eqref{uraba2}, and eventually to show that $\alpha=\pm 1,$ $\beta =0.$

\bl \la{l4}  Let $Q=F\circ T_s\circ (\alpha x+\beta)$, where $s\geq 5,$ $\alpha, \beta \in \R$, $\alpha\neq 0.$
\vskip 0.2cm
\noindent 1) If the coefficients $Ch_2(Q)$, $Ch_3(Q)$ vanish, then either $\alpha=\pm 1,$ $\beta=0,$ or 
\be \la{azxx} 4\beta^2=\frac{6}{(n-1)(2n-1)}, \ \ \ \alpha^2=\frac{2n-4}{2n-1}. \ee
\vskip 0.2cm
\noindent 2) If the coefficients $Ch_2(Q)$, $Ch_4(Q),$  vanish, then either $\alpha=\pm 1,$ $\beta=0,$ or 
\be \la{azxx+} 
4{\beta}^2=\frac{12}{(n-1)(2n-1)}, \ \ \ {\alpha}^2=\frac{2n-7}{2n-1}. \ee
In particular, in both cases $\alpha^2< 1$ and $\beta \neq 0,$ unless $\alpha=\pm 1,$ $\beta=0.$
\el

\pr Set $n=\deg Q.$ 
By Lemma \ref{l1}, for some $b_0\in \R$ we have: $$C_i(Q)=C_i(F\circ T_s\circ (\alpha x+\beta))=C_i(b_0T_{n/s}\circ T_s\circ (\alpha x+\beta))=C_i(b_0T_n\circ (\alpha x+\beta)), \ \ \ 0\leq i \leq s-1,$$ 
implying that 
$$Q=(b_0T_n)\circ (\alpha x+\beta)+R_2,$$ where $R_2$ ia a polynomial such that $\deg R_2\leq n-s$. 
Thus, if $Ch_2(Q)=0,$ we have:
\be \la{0} Q=c_0T_n+c_1T_{n-1}+c_3T_{n-3}+c_4T_{n-4}+ R_{1}=b_0T_n\circ (\alpha x+\beta)+R_{2},\ee
where  $\deg R_1, \deg R_2\leq n-5$, and $c_0,c_1,c_3,c_4,b_0\in \R$.
Changing $x$ to $x/2$, and $R_i(x)$ to $2R_i(x/2)$, we obtain a similar equality 
\be \la{4}  Q=c_0T_n^*+ c_1T_{n-1}^*+ c_3T_{n-3}^*+ c_4T_{n-4}^*+R_{1}={ b_0}T_n^{*}\circ (\t\alpha x+\t\beta)+R_{2},\ee
where $\t\beta=2\beta,$ $\t\alpha=\alpha.$ Furthermore, without loss of generality we may assume that $c_0=1,$ implying 
that $ b_0=1/{\t \alpha}^n$ and $c_1=\t\beta n/\t\alpha$. Thus, we can rewrite \eqref{4} in the form  
\be \la{5} Q=T_n^*+\frac{\t \beta n}{\t \alpha}T_{n-1}^*+c_3T_{n-3}^*+ c_4T_{n-4}^*+ R_{1}=\frac{1}{\t \alpha^n}T_n^{*}\circ (\t\alpha x+\t\beta)+ R_{2}.\ee

Calculating $C_2(Q),$ $C_3(Q),$ $C_4(Q)$ for both representations of $Q$ in \eqref{5} using
formula \eqref{f} (and the Taylor formula, for the second representation), we obtain the equalities 
$$-n=\frac{1}{(n-2)!{\t \alpha}^2}\left[\frac{n!{\t \beta^2}}{2!}-n(n-2)!\right],$$
$$\frac{-{\t \beta}n(n-1)}{{\t \alpha}}+ c_3=\frac{1}{(n-3)!{\t\alpha^3}}\left[\frac{n!{\t \beta^3}}{3!}-n(n-2)!\t \beta\right],$$ 
$$\frac{n(n-3)}{2}+c_4=\frac{1}{(n-4)!{\t\alpha}^4}\left[
\frac{n!{\t \beta}^4}{4!}-\frac{n(n-2)!{\t \beta}^2}{2!}+\frac{n(n-3)(n-4)!}{2}\right].$$

It follows from the first of these equalities that 
\be \la{gopa} {\t \alpha}^2=1-\frac{(n-1){\t \beta}^2}{2}.\ee If $c_3=0$, then 
substituting this value of ${\t \alpha}^2$ into third equation we   
obtain that either $\beta=0$ and then $\alpha=\pm 1$ by \eqref{gopa}, or \eqref{azxx} holds.
Similarly, if $c_4=0$, then 
substituting the value of ${\t \alpha}^2$ from \eqref{gopa} into the fourth equation we   
obtain 
$$\frac{n!{\t \beta}^4}{4!}-\frac{n(n-2)!{\t \beta}^2}{2!}+\frac{n(n-3)(n-4)!}{2}= \frac{(n-4)!n(n-3)}{2}\left[\frac{(n-1)^2{\t \beta}^4}{4}-(n-1){\t \beta}^2+1\right]
$$ implying that, either $\t \alpha=\pm 1,$ $\t \beta=0,$  or  \eqref{azxx+} holds.

Finally, clearly $\alpha^2< 1$, $\beta \neq 0,$ unless $\alpha=\pm 1,$ $\beta=0.$
\qed

\bl \la{dfg3} Let $Q\in C(n,m_1,m_2)$, where $m_1$ and $m_2$ are 
coprime, $m_1\geq 3,$ $m_2\geq 3$. Then at least one of the coefficients $Ch_2(Q),$ $Ch_4(Q),$ $Ch_6(Q)$ vanishes.
\el 
\pr Assume that $Ch_2(Q)\neq 0,$ $Ch_4(Q)\neq 0,$ and show that this implies that $Ch_6(Q)=0.$ First, $Ch_2(Q)\neq 0$ implies by \eqref{inv} that $m_2\mid n-2.$
It follows now from $Ch_4(Q)\neq 0$ by \eqref{inv} that $m_1=4.$
Therefore, $Ch_6(Q)=0$ for otherwise \eqref{inv} implies that $m_2=4$ in contradiction with $\GCD(m_1,m_2)=1.$ \qed

\bl \la{ll3} Let $Q=U\circ T_s\circ \alpha x$, where $U\in \R[x],$ $\alpha \in \R\setminus\{0\},$ and $s\geq 6.$
\vskip 0.2cm
\noindent 1) If the coefficient $Ch_2(Q)$ vanishes, then $\alpha^2=1.$
\vskip 0.2cm
\noindent 2) If the coefficient $Ch_4(Q)$ vanishes, then
either $\alpha^2=1,$ or $\alpha^2=\frac{n-3}{n-1}.$
\vskip 0.2cm
\noindent 3) If the coefficient $Ch_6(Q)$ vanishes, then $\alpha^2$ is a root of the equation 
\be \la{kon} (n^2-3n+2)t^2+(-2n^2+12n-16)t+(n^2-9n+20)=0.\ee In particular, in all the cases the inequality $\alpha^2< 1$ holds, unless $\alpha^2= 1$.
\el

\pr Set $n=\deg Q.$ As in Lemma \ref{l4} we arrive to the equality  
\be \la{00} T_n^{*}+c_2T_{n-2}^{*}+c_4T_{n-4}^{*} + c_6T_{n-6}^{*}+ R_1=\frac{1}{\t \alpha^n}T_n^{*}\circ (\alpha x)+ R_2,\ee
where $\deg R_1, \deg R_2\leq n-7$ (since obviously $c_1=c_3=c_5=0$). 
Calculating now $C_2(Q),$ $C_4(Q),$ $C_6(Q)$ for both representations of $Q$ in \eqref{00}, we obtain:
$$-n+ c_2=-\frac{n}{\alpha^{2}},$$
$$\frac{n(n-3)}{2}-(n-2)c_2+c_4=\frac{n(n-3)}{2\alpha^4},$$
$$-\frac{n(n-4)(n-5)}{6}+\frac{(n-2)(n-5)c_2}{2}-(n-4)c_4+c_6=-\frac{n(n-4)(n-5)}{6\alpha^6}$$

It follows from the first equality that $$c_2=\frac{n(\alpha^2-1)}{\alpha^2},$$ implying that
if $c_2=0$, then $\alpha^2=1.$ Substituting now the value of $c_2$ into the second equality 
we obtain that 
$$c_4=\frac{n(n-3)}{2\alpha^4}-\frac{n(n-3)}{2}+\frac{n(n-2)(\alpha^2-1)}{\alpha^2}=n\left(\frac{(n-1)\alpha^4-2(n-2)\alpha^2+(n-3)}{2\alpha^4}\right).
$$ Since $$(n-1)\alpha^4-2(n-2)\alpha^2+(n-3)=(n-1)(\alpha^2-1)\left(\alpha^2-\frac{n-3}{n-1}\right),$$
this implies that if $c_4=0$, then either $\alpha^2=1,$ or $\alpha^2=\frac{n-3}{n-1}.$

Finally, $$c_6=-\frac{n(n-4)(n-5)}{6\alpha^6}+\frac{n(n-4)(n-5)}{6}-\frac{(n-2)(n-5)c_2}{2}+(n-4)c_4=$$
$$=-\frac{n(n-4)(n-5)}{6\alpha^6}+\frac{n(n-4)(n-5)}{6}-\frac{(n-2)(n-5)}{2}\frac{n(\alpha^2-1)}{\alpha^2}+$$ 
$$+(n-4)n\left(\frac{(n-1)\alpha^4-2(n-2)\alpha^2+(n-3)}{2\alpha^4}\right)=$$
$$
=n\left(\frac{(n^2-3n+2)\alpha^6+(-3n^2+15n-18)\alpha^4+(3n^2-21n+36)\alpha^2-n^2+9n-20}{6\alpha^6}\right).
$$
Now the statement follows from the factorization  
$$
(n^2-3n+2)\alpha^6+(-3n^2+15n-18)\alpha^4+(3n^2-21n+36)\alpha^2-n^2+9n-20=$$ $$=(\alpha^2-1)\left((n^2-3n+2)\alpha^4+(-2n^2+12n-16)\alpha^2+(n^2-9n+20)\right). $$

In order to finish the proof we only must show that the absolute values of roots of  
equation \eqref{kon} are less than one. Solving \eqref{kon}, we find two roots 
$$t_1=\frac{n^2-6n+8-\sqrt{3n^2-18n+24}}{n^2-3n+2}=\frac{\sqrt{(n-2)(n-4)}(\sqrt{(n-2)(n-4)}-\sqrt{3})}{(n-1)(n-2)},
$$
$$t_2= \frac{n^2-6n+8+\sqrt{3n^2-18n+24}}{n^2-3n+2}=\frac{\sqrt{(n-2)(n-4)}(\sqrt{(n-2)(n-4)}+\sqrt{3})}{(n-1)(n-2)}.$$
Clearly, for $n\geq 6$ the inequality $0<t_1<t_2$ holds. Finally, 
$$t_2-1=\frac{\sqrt{3(n-2)(n-4)}-3n+6}{(n-1)(n-2)}=\frac{\sqrt{3(n-2)}(\sqrt{n-4}-\sqrt{3(n-2)})}{(n-1)(n-2)}<0. \ \ \ \Box$$

\vskip 0.2cm
\noindent{\it Proof of Proposition  \ref{vot}.}
Assume first that $Ch_1(Q)\neq 0.$
Then Lemma \ref{dfg} implies that either $Ch_2(Q)$, $Ch_3(Q)$ vanish, or $Ch_2(Q)$, $Ch_4(Q)$ vanish.
It follows now from Lemma \ref{l4} that, unless $\alpha=\pm 1,$ $\beta =0,$ the conditions $\alpha<1$, $\beta\neq 0$ hold. So, assume that $\alpha<1$, $\beta\neq 0$.

Rewrite equality \eqref{uraba2} in the form  
\be \la{uraba3} P=\t V_1\circ T_{\t m_1}+\t V_2\circ T_{\t m_2}=U\circ T_{m_1m_2}\circ \left(\frac{x-\beta}{\alpha}\right).\ee Since $\beta \neq 0,$ Corollary \ref{c1} 
applied to \eqref{uraba3} implies that $C_1(P)\neq 0.$ Applying now Lemma \ref{dfg} and Lemma \ref{l4} to \eqref{uraba3} in the same way as before to equality \eqref{uraba}, we  conclude that 
$1/\alpha<1.$ The contradiction obtained proves that $\alpha=\pm 1,$ $\beta =0.$

Assume now that $Ch_1(Q)= 0.$ Then $\beta=0$, by Corollary \ref{c1}.  
Furthermore, by Lemma \ref{dfg3}
at least one of the coefficients $Ch_2(Q),$ $Ch_4(Q),$ $Ch_6(Q)$ vanish,
implying by Lemma \ref{ll3} that, unless $\alpha=\pm 1,$ $\beta =0,$ the condition $\alpha<1$ holds. Since $\beta=0$ implies by Corollary \ref{c1} that $C_1(P)=0$ in view of \eqref{uraba3}, the assumption $\alpha<1$  leads to a contradiction in the same way as before. \qed

\subsection{Proof of Theorem \ref{t2} in the cases (3,1), (3,2), (3,3).}
The  case (3,1) reduces  to the case (2,1) as follows. 
We start from the equality 
\be \la{ura} Q=V_1\circ T_{2m_1}+V_2\circ T_{2m_2}+V_3\circ xR(x^2)\circ T_{m_1m_2}
=\t U\circ 
 x^2{\t R}^2(x^2)\circ(\alpha x+\beta),\ee where $V_1,V_2,V_3,R,\t R, \t U\in \R[x]$, $\alpha,\beta\in \R,$ $\alpha\neq 0,$
and $m_1\geq 3,$ $m_2\geq 3$ are coprime and odd. 
It follows from the first representation for $Q$ in \eqref{ura} that $Q$ can be written in the form 
\be \la{ebi} Q=d_nT_n+d_{n-1}T_{n-1}+\dots +d_{1}T_1+d_0,\ \ \ d_i\in \R,\ee
where $d_i=0$, unless $i$ is divisible either by $2m_1,$ or by $2m_2$, or by $m_1m_2$.
Clearly, conditions imposed on $m_1,m_2$ imply that 
\be \la{dif1} \vert 2m_1k_1-2m_2k_2\vert\geq 2, \ \ \ {\rm unless} \ \ \ \vert 2m_1k_1-2m_2k_2\vert=0,\ee  
\be \la{dif1+} \vert 2m_ik_1-2m_ik_2\vert\geq  2m_i\geq  6, \ \ \ {\rm unless} \ \ \ \vert 2m_ik_1-2m_ik_2\vert=0,\ \ \ i=1,2\ee  
\be \la{dif2} \vert 2m_ik_1-m_1m_2k_2\vert\geq  m_i\geq 3, \ \ \ {\rm unless} \ \ \ \vert 2m_ik_1-m_1m_2k_2\vert=0,\ \ \ i=1,2\ee
\be \la{dif3} \vert m_1m_2k_1-m_1m_2k_2\vert\geq  m_1m_2 \geq 15, \ \ \ {\rm unless} \ \ \ \vert m_1m_2k_1-m_1m_2k_2\vert=0.\ee
Therefore, $Ch_1(Q)=0$, implying that $C_1(Q)=0$, since $C_1(T_n)=0.$ 
It follows now from the second representation for $Q$ in \eqref{ura} by Corollary \ref{c1} that $\beta=0$.
Since the polynomial $\t W_1=z^2\circ (\alpha x+\beta)$ satisfies $\t W_1(a)=\t W_1(b)$, this implies that 
$a=- b.$ Therefore, solution $P,Q$ also belongs to the case (2,1) considered earlier (see the remarks after Theorem \ref{mpc}).

In the case (3,2) 
there exist  $V_1,V_2,V_3, U, R, \t V_1, \t V_2, \t U\in \R[x]$ and $\alpha,\beta\in \R$, $\alpha\neq 0,$  such that 
\be \la{ura+} Q=V_1\circ T_{2m_1}+V_2\circ T_{2m_2}+V_3\circ xR(x^2)\circ T_{m_1m_2}
=\t U\circ T_{\t m_1\t m_2} \circ (\alpha x+\beta),\ee 
\be \la{rau+} P=\t V_1\circ T_{\t m_1} \circ  (\alpha x+\beta)+\t V_2\circ T_{\t m_2}\circ  (\alpha x+\beta)=U\circ x^2R^2(x^2)\circ T_{m_1m_2}\ee
where $m_1\geq 3$, $m_2\geq 3$ are odd, $\GCD(m_1,m_2)=1,$ and $\t m_1\geq 2$, $\t m_2\geq 2$, $\GCD(\t m_1,\t m_2)=1,$
Further, without loss of generality we may assume that $a\neq - b$, for otherwise $P,Q$ belongs to the case $(2,2).$ Besides, we may assume that $\t m_1\geq 3$, $\t m_2\geq 3,$ for otherwise $P,Q$ belongs to the case $(3,1).$

Since equalities \eqref{ura+}, \eqref{rau+} may be written in the form  
$$ Q=(V_1\circ T_2+ V_3\circ xR(x^2)\circ T_{m_2})\circ T_{m_1} +(V_2\circ T_2)\circ T_{m_2}
=\t U\circ T_{\t m_1\t m_2} \circ (\alpha x+\beta),$$ 
$$ P=\t V_1\circ T_{\t m_1} \circ  (\alpha x+\beta)+\t V_2\circ T_{\t m_2}\circ  (\alpha x+\beta)=(U\circ x^2R^2(x^2))\circ T_{m_1m_2},$$ 
it follows from Proposition \ref{vot} that $\alpha=\pm 1,$ $\beta =0.$ 
Since we assumed that $a\neq - b,$ it follows now from the equalities 
$$T_{2m_1}(a)=T_{2m_1}(b), \ \ \ T_{2m_2}(a)=T_{2m_2}(b), \ \ \ T_{{\t m}_1}(\pm a)=T_{{\t m}_1}(\pm b)$$ 
by  Lemma \ref{skun}, b), taking into account \eqref{kots}, that 
$T_s(a)=T_s(b)$, where either $s=\GCD(2m_1, {\t m}_1)$, or $s=\GCD(2m_2, {\t m}_2).$ Finally, 
since
$$Q=\t U\circ T_{\t m\t n}\circ (\pm z)=\t U\circ (\pm  T_{\t m\t n})=\t U\circ (\pm  T_{\t m\t n/s})\circ T_s$$
and 
$$P=U\circ x^{2}R^2(x^{2})\circ  T_{m_1m_2}=U\circ x^{2}R^2(x^{2})\circ  T_{m_1m_2/s}\circ T_s,$$
we conclude that \eqref{c} holds for $W=T_s.$

The proof in the case (3,3) is similar: 
there exist  $V_1,V_2,V_3, U, R, \t V_1, \t V_2, \t V_3, \t U, \t R\in \R[x]$ and $\alpha,\beta\in \R$, $\alpha\neq 0,$  such that 
$$ Q=V_1\circ T_{2m_1}+V_2\circ T_{2m_2}+V_3\circ xR(x^2)\circ T_{m_1m_2}
=\t U\circ x^2\t R^2(x^2)\circ T_{\t m_1\t m_2} \circ (\alpha x+\beta),$$ 
$$ P=\t V_1\circ T_{\t 2m_1} \circ  (\alpha x+\beta)+\t V_2\circ T_{\t 2m_2}\circ  (\alpha x+\beta)+\t V_3\circ x\t R(x^2)\circ T_{\t m_1\t m_2}
=U\circ x^2R^2(x^2)\circ T_{m_1m_2}$$
where $m_1\geq 3$, $m_2\geq 3$, $\t m_1\geq 3$, $\t m_2\geq 3$ are odd, $\GCD(m_1,m_2)=1,$ $\GCD(\t m_1,\t m_2)=1.$ Moreover,
without loss of generality we may assume that $a\neq - b$. 
Further, using  Proposition \ref{vot} we conclude as above that  $\alpha=\pm 1,$ $\beta =0.$  Finally, it follows from 
the equalities 
$$T_{2m_1}(a)=T_{2m_1}(b), \ \ \ T_{2m_2}(a)=T_{2m_2}(b), \ \ \ T_{{\t 2m}_1}(\pm a)=T_{{\t 2m}_1}(\pm b)$$ 
that \eqref{c} holds for $W=T_s,$  where either $s=\GCD(2m_1, {2\t m}_1)$, or $s=\GCD(2m_2, 2{\t m}_2).$


\end{document}